\documentclass{article} 
\usepackage[dvips]{graphicx}
\usepackage{fancyhdr}
\usepackage{amsmath,amssymb,latexsym}
\usepackage{color}

\newcommand{\norm}[2]{
\left\| #2 \right\|_{#1}
}
\newcommand{\Hil}[0]{
\mathcal{H} 
}

\newcommand{\CC}[0]{
\mathbb{C} 
}

\newtheorem{def.}[theorem]{Definition}
\newtheorem{prop.}[theorem]{Proposition}
\newtheorem{lem.}[theorem]{Lemma}
\newtheorem{cor.}[theorem]{Corollary}
\newtheorem{conj.}[theorem]{Conjecture}
\newtheorem{com.}[theorem]{Commentary}
\newtheorem{Bsp.}{Example}[section]
\newtheorem{rem.}[theorem]{Remark}
\newtheorem{rems.}[theorem]{Remarks}

\newenvironment{proof}{\noindent \bf Proof: \rm}{$ \hspace{\stretch{1}} \Box $ \vspace{5mm}}

\def\dom{{\rm dom}}
\def\ran{{\rm ran}}
\def\span{{\rm span}}
\def\<{\left<}
\def\>{\right>}
\newcommand{\mn}{\mathbb N}

\begin{document}

\title{\bf\vspace{-39pt} 
Classification of General Sequences by Frame-Related Operators}

\author{Peter Balazs \\ \small $^{a)}$ Acoustics Research Institute, \\
\small Wohllebengasse 12-14,  Vienna A-1040, Austria \\ \small Peter.Balazs@oeaw.ac.at\\
\\  
Diana T. Stoeva$^{a),\, b)}$ \\ \small $^{b)}$ Department of Mathematics, \\ 
\small University of Architecture, Civil Engineering and Geodesy, \\ \small Blvd Hristo Smirnenski 1, 1046 Sofia, Bulgaria \\ \small  std73std@yahoo.com \\
\\
Jean-Pierre Antoine \\ \small Institut de Recherche en Math\'ematique et  Physique, \\ 
\small Universit\'e catholique de Louvain,   \\
\small B - 1348   Louvain-la-Neuve, Belgium  \\ \small Jean-Pierre.Antoine@uclouvain.be }

\maketitle

\begin{abstract}
This note is a survey and collection of results, as well as presenting some original research.
For Bessel sequences and frames, the analysis, synthesis and frame operators as well as the Gram matrix
are well-known, bounded operators. 
We investigate these operators for arbitrary sequences, which in general lead to possibly unbounded operators. 
We characterize various classes of sequences in terms of these operators and vice-versa. 
Finally, we classify these sequences by operators applied on orthonormal bases.
\vspace{5mm}  
\vspace{3mm}
\end{abstract}

\section{Introduction}

Frames \cite{duffschaef1, Casaz1,ole1} 
have become an important topic in applied mathematics in the last decades. They are also used in many engineering applications, for example in signal processing \cite{boelc1}.
Clearly, they have some important advantages compared to orthonormal bases. 
However, 
even the frame condition cannot always be satisfied for the whole space, and so other classes of sequences have been investigated, for example, frame sequences, Bessel sequences, 
lower frame sequences, 
and Riesz-Fischer sequences  \cite{ xxlelgeb07, casoleli1, olepinv, ole1}.
For such sequences, which need not be frames in general, 
the frame-related operators, i.e. the analysis, the synthesis and the frame operator can still be defined, see e.g. \cite{casoleli1,olepinv}. 
In general, these operators can be unbounded.

In this note we want to give an overview of the connection between the properties of those operators and those of the sequences.
We collect some existing results, extend them and add new, original results. 
Some known results are proved here for the sake of completeness, in  particular, when there is no proof in the literature. 
Note that one result was obtained independently in \cite{elgebunb1}. While the paper \cite{elgebunb1} focuses on the investigation of sufficient conditions for the closability of the synthesis operator,
 the main aim of the present paper is to give a general overview of all the associated operators.

Our notation and some preliminary results are given in Section \ref{sec:prelnot0}. 
In Section \ref{sec:unbframop0} we consider the operators associated to an arbitrary sequence, namely, the analysis, 
synthesis and \lq frame\rq \ operators, as well as the operator based on the Gram matrix, and we investigate 
their properties.
Section \ref{sec:framclass0} is devoted to 
the characterization of sequences (Bessel and frame sequences, frames, Riesz bases, lower frame sequences, Riesz-Fischer sequences, complete sequences) 
 via the associated operators. 
Finally, we consider operators that preserve the sequence type and present the classification via orthonormal bases.

\section{Notation and Preliminaries} \label{sec:prelnot0}

Throughout the paper we consider a sequence $\Psi=(\psi_k)_{k=1}^\infty$ with elements from a (infinite dimensional) Hilbert space $(\Hil, \<\cdot,\cdot\>)$. 
The notation $\span{\{\Psi\}}$ is used to denote the linear span of $(\psi_k)_{k=1}^\infty$ and $\overline{\span}{\{\Psi\}}$ denotes the closed linear span of $(\psi_k)_{k=1}^\infty$. 
The sequence $(e_k)_{k=1}^\infty$ denotes an orthonormal basis of $\Hil$ and $(\delta_k)_{k=1}^\infty$ denotes the canonical basis of $\ell^2$. 
The notion {\it operator} is used for a linear mapping. 
Given an operator $F$, 
we denote its domain by $\dom(F)$, its range by $\ran(F)$ and its kernel by $\ker(F)$.
First we recall the definitions of the basic concepts used in the paper.

\begin{def.} The sequence $\Psi$ is called
\begin{itemize}

\item  complete in $\Hil$ if  $\overline{\span}{\{\Psi\}}=\Hil$;

\item 
 a Bessel sequence for $\Hil$ with bound $B$ if  $B>0$ and   
 $\sum_{k=1}^\infty |\<f,\psi_k\>|^2 \leq B\|f\|^2$ for every $f\in\Hil$; 
 
\item a lower frame sequence for $\Hil$ with bound $A$ if $A>0$ and $A\|f\|^2 \leq \sum_{k=1}^\infty |\<f,\psi_k\>|^2$ for every $f\in\Hil$;

\item  a frame for $\Hil$ with bounds $A,B$ if it is a Bessel sequence for $\Hil$ with bound $B$ and a lower frame sequence for $\Hil$ with bound $A$;

\item  a frame-sequence if it is a frame for its closed linear span;

\item  a Riesz basis for $\Hil$ with bounds $A,B$ if $\Psi$ is complete in $\Hil$,
$A>0$, $B<\infty$, and $A\sum |c_k|^2 \leq \|\sum c_k\psi_k\|^2\leq B\sum |c_k|^2$ for all finite scalar sequences $(c_k)$ (and hence, for all $(c_k)_{k=1}^\infty\in\ell^2$);

\item  a Riesz-Fischer sequence with bound $A$ if $A>0$ and $A\sum |c_k|^2 \leq \|\sum c_k\psi_k\|^2$ for all finite scalar sequences $(c_k)$ (and hence, for all $(c_k)_{k=1}^\infty\in\ell^2$ such that $\sum_{k=1}^\infty c_k\psi_k$ converges in $\Hil$).

\end{itemize}

\end{def.}

\begin{lem.} \label{sec:bansteinhltwo1} Let $(d_k)_{k=1}^\infty$ be a sequence such that 
\begin{equation} \label{conv1}
\sum \limits_{k=1}^{\infty} c_k d_k \mbox{ converges  for all } (c_k)_{k=1}^\infty \in l^2.
\end{equation}
 Then $(d_k)_{k=1}^\infty \in l^2$. 
\end{lem.}
\begin{proof} If (\ref{conv1}) is assumed to hold with absolute convergence of the series, the above statement is proved in \cite[\S 30 1.(5)]{koeth1}. 
 Similar proof (using $x_j=\frac{1}{i|M_i|}\overline{v_j}$ instead of $x_j=\frac{1}{iM_i}v_j$) leads to the validity of the above lemma.
\end{proof}

We also need the following well-known results: 
\begin{prop.} \label{sec:closeunb1} Let an operator $F$ be densely defined. Then 
\begin{itemize} 

\item[{\rm (a)}] 
 {\rm  \cite[Prop.\,X.1.6]{conw1}} 
 $F^*$ is closed; in particular, every self-adjoint operator is closed;
 
\item[{\rm (b)}] 
{\rm  \cite[Prop.\,X.1.6]{conw1}} $F$ is closable if and only if $F^*$ is densely defined;

\item[{\rm (c)}]
{\rm  \cite[Th.\,VIII.1]{rs1}} If $F$ is closable, then its closure is $\overline{F} = (F^{*})^{*}$ and $\overline{F}^{\, *} = F^*$;

\item[{\rm (d)}] 
{\rm  \cite[Prop. X.1.13]{conw1}} ${\left( \ran(F) \right)}^\bot= \ker (F^*)$; 

\item[{\rm (e)}] 
{\rm  \cite[Prop. X.4.2]{conw1}} 
If $F$ is closed, then $F^*F$ is self-adjoint.
\end{itemize}
\end{prop.}

\section{Associated Operators for Arbitrary Sequences} \label{sec:unbframop0}

For any sequence $\Psi$, 
the associated analysis operator $C$, synthesis operator $D$, \lq frame\rq \ operator $S$, and Gram operator $G$
are (possibly unbounded) operators, defined as  follows:

\vspace{.05in}
$C : \dom(C)\to \ell^2$, $ C f = 
\left( \left< f , \psi_k \right> \right)_{k=1}^\infty$,  where 

\vspace{.05in}
\hspace{.2in} $ \dom(C) = \left\{ f \in \Hil \ : \ 
\left( \left< f , \psi_k \right> \right)_{k=1}^\infty \in\ell^2 \right\};$

\vspace{.1in}
$D : \dom(D)\to \Hil$, $ D (c_k)_{k=1}^\infty = 
\sum_{k=1}^\infty c_k \, \psi_k $, where 

\vspace{.05in}
\hspace{.2in}  $\dom(D) = \left\{ (c_k)_{k=1}^\infty \in \ell^2 \ : \ \sum_{k=1}^\infty c_k \, \psi_k \mbox{  converges in $\Hil$} \right\}$;

\vspace{.1in}
$S : \dom(S)\to \Hil$, $S f = \sum_{k=1}^\infty \left< f , \psi_k\right> \psi_k $,  where  

\vspace{.05in}
\hspace{.2in} $\dom(S) = \left\{f\in\Hil \ : \ \sum_{k=1}^\infty \left< f , \psi_k\right> \psi_k  \mbox{ converges in $\Hil$} \right\};$

\vspace{.1in}
$G : \dom(G)\to \ell^2$, $ G (c_k)_{k=1}^\infty = \left(\sum_{l=1}^\infty G_{k,l} c_l\right)_{k=1}^\infty$, where 

\vspace{.05in}
\hspace{.2in} $\dom(G) = \left\{ (c_k)_{k=1}^\infty \in \ell^2 \, : \, \sum_{l=1}^\infty G_{k,l} \, c_l \mbox{ converges } \forall k\in\mn \right.$

  \hspace{2.5in} 
$\left. \mbox{ and }  \left(\sum_{l=1}^\infty G_{k,l} \, c_l\right)_{k=1}^\infty \in \ell^2 \right\},$
 
\vspace{.05in}   \hspace{.2in} and the Gram matrix $(G_{k,l})_{k,l}$ is defined by $G_{k,l} = \left< \psi_l, \psi_k \right>$, $k,l\in\mn$.

\subsection{Domains}

\begin{lem.} \label{dom1} Given an arbitrary sequence $\Psi$, the following statements hold.
\begin{itemize}
\item[{\rm (i.1)}] $\dom(S) \subseteq \dom(C)$.
\item[{\rm (i.2)}] $\dom(S)=\dom (C)$ if and only if  $\ran(C)\subseteq \dom(D)$. 

\item[{\rm (ii)}] $\dom(D)\subseteq \dom(G)$ if and only if $\ran(D)\subseteq \dom(C)$.

\item[{\rm (iii)}] $D$ is densely defined. 

\item[{\rm (iv)}] {\rm \cite[Prop. 4.5]{olepinv}} If $\sum_{l=1}^\infty |\<\psi_l,\psi_k\>|^2 <\infty$ for every $k\in\mn$, then $C$ is densely defined. 

\item[{\rm (v)}] If $\sum_{l=1}^\infty |\<\psi_l,\psi_k\>|^2 <\infty$ for every $k\in\mn$, then $G$ is densely defined.

\item[{\rm (vi)}] 
 If $\sum_{l=1}^\infty \sum_{k=1}^\infty |\<\psi_l,\psi_k\>|^2 <\infty$, then $\dom(G) = \dom(D)=\ell^2$, $G$ is a Hilbert-Schmidt operator and $\Psi$ is a Bessel sequence for $\Hil$. 

\item[{\rm (vi)}] 
 If $\sum_{l=1}^\infty \sum_{k=1}^\infty |\<\psi_l,\psi_k\>|^2 <\infty$, then $\dom(G) = \dom(D)=\ell^2$ and $G$ is a Hilbert-Schmidt operator. 
 
\item[{\rm (vii)}] If $\psi_k\in\dom(S)$ for every $k\in\mn$, then $S$ is densely defined.

\end{itemize}
\end{lem.}

\begin{proof} (i.1)
Assume that $f \in \dom(S)$. Then $\sum_{k=1}^N \<\left< f , \psi_k\right> \psi_k,f\> \to \<Sf,f\>$ as $N\to\infty$, which implies that 
$\sum_{k=1}^N \left| \left< f , \psi_k\right> \right|^2$ converges as $N\to\infty$ and thus $f\in\dom(C)$.

(i.2) If $\ran(C)\subseteq \dom(D)$, it is obvious that $\dom(C)\subseteq \dom (S)$ and the other inclusion is given in (i.1). 
Now assume that $\dom(S)=\dom (C)$ and take $Cf\in \ran(C)$. Then $f\in\dom(S)$, which implies that $Cf\in\dom(D)$.

(ii) First observe that if $\sum_{l=1}^\infty c_l \psi_l$ converges in $\Hil$, then $\sum_{l=1}^\infty c_l \<\psi_l,\psi_k\>$ converges for every $k\in\mn$ and $\sum_{l=1}^\infty c_l \<\psi_l,\psi_k\>=\<Dc,\psi_k\>$, $\forall k\in\mn$.

Assume that $\dom(D)\subseteq \dom(G)$ and take $Dc\in\ran(D)$. Then $c\in\dom(G)$, which implies 
that $(\<Dc,\psi_k\>)_{k=1}^\infty \in\ell^2$ and thus, $Dc\in\dom(C)$.

Now assume that $\ran(D)\subseteq \dom(C)$ and take $c\in\dom(D)$. Then $Dc\in\dom(C)$ and 
 $\sum_{l=1}^\infty c_l \<\psi_l,\psi_k\>=\<Dc,\psi_k\>$, $\forall k\in\mn$. Therefore, $c\in\dom(G)$.

(iii) is clear, since 
 the finite sequences are dense in $\ell^2$.

(iv) Since the statement is given in \cite{olepinv} without proof, for the sake of completeness we include a proof here. First observe that $(\overline{\span}\{\Psi\})^{\bot}\subseteq \dom(C)$. Assume that $\sum_{l=1}^\infty |\<\psi_l,\psi_k\>|^2 <\infty$ for every $k\in\mn$. Then $\span\{\Psi\}\subseteq \dom(C)$. Therefore, $\Hil=\overline{\span}\{\Psi\} \oplus (\overline{\span}\{\Psi\})^{\bot} \subseteq \overline{\dom(C)}$, which implies that $\dom(C)$ is dense in $\Hil$.

(v) 
Let $\sum_{l=1}^\infty |\<\psi_l,\psi_k\>|^2 <\infty$ for every $k\in\mn$. 
Then all the canonical vectors $\delta_k$ (and thus, all the finite sequences) belong to $\dom(G)$. Then the conclusion follows as in (iii).

(vi) 
Assume that $\sum_{l=1}^\infty \sum_{k=1}^\infty |\<\psi_l,\psi_k\>|^2 <\infty$. 
Then $G$ is a 
Hilbert-Schmidt operator \cite[Cor. 16.11]{MeiseVogt}.
Furthermore, $\Psi$ is a Bessel sequence \cite[Lemma 3.5.1]{ole1} (see Prop. \ref{sec:classgram1}(a2)), which implies that $\dom(D)=\ell^2$  \rm\cite[Theor. 3.2.3]{ole1} (see Prop. \ref{sec:classwsynth1}).

(vii) Follows in a similar way as in (iv), having in mind that $(\overline{\span}\{\Psi\})^{\bot}\subseteq \dom(S)$. 
\end{proof}

\begin{rems.} 
\label{rem32} \

{\rm

\vspace{.05in}
(1) Concerning Lemma \ref{dom1}(i):

If $\Psi$ is a Bessel sequence for $\Hil$, then $\dom(S)=\dom(C)=\Hil$. Note that the equality $\dom(S)=\dom(C)$ might hold even in cases when $\Psi$ is not a Bessel sequence. Consider for example the non-Bessel sequence 

\centerline{$\Psi=(e_1, e_1, e_2, e_1, e_2, e_3, e_1, e_2, e_3, e_4, \ldots)$.}

\noindent Then $\dom(C)=\{0\}$ and thus, Lemma \ref{dom1}(i.1) implies that $\dom(S)=\{0\}$. 
\label{nullex}

 Note that the equality $\dom(S)= \dom(C)$ might fail. 
 Consider for example the non-Bessel sequence  $\Psi=(\frac{1}{2}e_1, 2e_2, \frac{1}{2^2}e_1, 2^2 e_3, \frac{1}{2^3}e_1, 2^3 e_4, \ldots)$. As an example of an element which belongs to $\dom(C)$ and does not belong to $\dom(S)$ consider  $h\in\Hil$ such that $\<h,e_n\>=2^{-2(n-1)}$,  
  $\forall n\in\mn$.

\vspace{.05in}
(2) Concerning Lemma \ref{dom1}(ii): 

If $\Psi$ is a Bessel sequence for $\Hil$, then $\dom(D)=\dom(G)=\ell^2$. 
For a non-Bessel sequence $\Psi$, the inclusion $\dom(D)\subseteq \dom(G)$ might fail.
Take for example $\Psi=(e_1, e_2, e_1, e_3, e_1, e_4, \ldots)$. Then 
$\delta_1\in\dom(D)$ and $\delta_1\notin\dom(G)$.

Note that a necessary condition for the validity of $\dom(G)=\dom(D)$ is that all the canonical vectors are in $\dom(G)$, which holds if and only if $\sum_{l=1}^\infty |\<\psi_l,\psi_k\>|^2<\infty$ for every $k$.

\vspace{.05in}
(3) Concerning Lemma \ref{dom1}(iii): 

If 
$\dom(D)=\ell^2$, then 
$\sum_{k=1}^\infty c_k\psi_k$ converges unconditionally for every $(c_k)\in\ell^2$ 
and $\Psi$ is a Bessel sequence for $\Hil$
(see \cite[Cor. 3.2.4 and 3.2.5]{ole1}).  
Now one could wonder whether this is pointwise true for an arbitrary sequence, i.e. whether
$c \in \dom(D)$ implies that $\sum_{k=1}^\infty c_k\psi_k$ is already unconditionally convergent. This 
is not true in general, as shown by the following easy example. Consider the non-Bessel sequence $\Psi=(e_1,e_1,e_1,\ldots)$ and $(c_k)=(1, -\frac{1}{2}, \frac{1}{3}, -\frac{1}{4}, \ldots)\in\ell^2$. Then $\sum_{k=1}^\infty c_k\psi_k = (\ln2) e_1$.   
However, if $(\epsilon_k)=(1,-1,1,-1,\ldots)$, then the series $\sum_{k=1}^\infty c_k \epsilon_k \psi_k= \sum_{k=1}^\infty \frac{1}{k} e_1$ does not 
converge, which by \cite[Theor. 2.8]{he98-1} implies that $\sum_{k=1}^\infty c_k\psi_k$ is conditionally convergent.

\vspace{.05in}
(4) Concerning Lemma \ref{dom1}(iv),(v): 

The assumption in (iv) and (v) is clearly satisfied when $\Psi$ is a Bessel sequence for $\Hil$ (in which case $\dom(C)=\Hil$ and $\dom(G)=\ell^2$). However, the validity of $\sum_{l=1}^\infty |\<\psi_l,\psi_k\>|^2 <\infty$ for every $k\in\mn$ does not require $\Psi$ to be a Bessel sequence for $\Hil$, consider for example $\Psi=(e_1, 2e_2, 3e_3, \ldots)$.

\vspace{.05in} 
(5) Concerning Lemma \ref{dom1}(vi): 

For frames, the assumption in (vi) is equivalent to the Hilbert space being finite dimensional \cite{xxlfinfram1}. 
 This is not necessarily valid any more for other sequences,  consider for example $\Psi=(e_1, e_2/2, e_3/3,e_4/4,...)$, it satisfies (vi) and $\Hil$ is infinite dimensional. 
} 
\end{rems.}

\subsection{Connection between the associated operators}

Having defined the operators 
associated to an arbitrary sequence, we will now investigate the relationships among them.

\begin{prop.} \label{connectop} For an arbitrary sequence $\Psi$, the following statements hold.  
\begin{itemize}
 \item[{\rm (i)}] $C=D^*$ \footnote{This result was obtained independently in \cite{elgebunb1}.}.  
 
 \item[{\rm (ii)}] {\rm \cite[Prop. 4.6]{olepinv}} If $C$ is densely defined, then $D\subseteq C^*$.
\item[{\rm (iii)}]  $S = D  C$. 
 \item[{\rm (iv)}] $  C D \subseteq G$.  
   \item[{\rm (v)}] {\rm \cite[Lemma 3.1]{casoleli1}} $C$ is closed. 
      \item[{\rm (vi)}] $S$ is symmetric\,\footnote{The terminology is not uniform in the literature. We use the definition of a symmetric operator given by 
      \cite[Def. XII.1.7]{DunfordS2}, \cite[Def. 13.3]{rudinengl},
      namely, $S$ is symmetric if $\<Sf,g\>=\<f,Sg\>, \forall f,g\in\dom(S)$, without the assumption of a dense domain.}, but not necessarily densely defined; it is positive on $dom(S)$ and positive definite on $\overline{\span} \{ \Psi\} \cap \dom(S)$.

\item[{\rm (vii)}] If $C$ is densely defined, then $S$ is closable.
      
  \item[{\rm (viii)}] $D$ is closable if and only if $C$ is densely defined.
   \item[{\rm (ix)}] $D$ is closed if and only if $C$ is densely defined and $D=C^*$.

   \item[{\rm (x)}] If $D$ is closed, then $G$ is densely defined. 
   
 \item[{\rm (xi)}]   If $\sum_{l=1}^\infty |\<\psi_l,\psi_k\>|^2 <\infty$ for every $k$, then $G$ is closed.
 \end{itemize}
\end{prop.}

\begin{proof} 
(i) 
First we show that $\dom(D^*) \subseteq \dom(C)$.
Fix $f \in \dom(D^*)$. This means that ${\mathfrak D} : \dom(D) \rightarrow \CC$, ${\mathfrak D} (c) := \left< D c , f\right>$, is a bounded functional.  
Since $\dom(D)$ is dense in $\ell^2$, there is a unique bounded extension ${\mathfrak D}_0 : \ell^2 \rightarrow \CC$. 
For $c=(c_k)_{k=1}^\infty \in \ell^2$, 
denote 
$$c_{_N} = (c_{_1}, c_{_2}, \dots, c_{_N}, 0, 0, \ldots ) \in c_{00} \subseteq \dom(D).$$
Clearly, $c_{_N} \rightarrow c$ in $\ell^2$-norm as $N\to\infty$. Hence, ${\mathfrak D}_0 c_{_N} \rightarrow {\mathfrak D}_0 c$ as $N\to\infty$. Therefore,  ${\mathfrak D}_0 c_{_N}=\sum_{k=1}^{N} c_k  \left< \psi_k , f \right>_{\Hil}$ converges as $N\to\infty$ for every $c\in\ell^2$. 
Now Lemma \ref{sec:bansteinhltwo1} implies that $\left( \left< \psi_k , f \right>_{\Hil} \right) \in \ell^2$, which proves that $f$ belongs to $\dom(C)$.

Now we show that $\dom(C)\subseteq \dom(D^*)$. 
Let $f \in \dom(C)$. For every $c=(c_k)_{k=1}^\infty \in \dom(D)$, 
\begin{equation} \label{cdstar} \left< Dc , f \right>_{\Hil} = 
\lim \limits_{N \rightarrow \infty} \left< \sum \limits_{k=1}^{N} c_k \, \psi_k , f \right>_{\Hil} = 
 \lim \limits_{N \rightarrow \infty}  \sum \limits_{k=1}^{N} c_k  \left< \psi_k , f \right>_{\Hil} 
 = \left< c , C f \right>_{\ell^2}.
 \end{equation}
Therefore, the functional $c 
\mapsto \left< Dc , f \right>_{\Hil}$ is bounded on $\dom(D)$ by $\|C f\|$. Hence, $f \in \dom(D^*)$.

Furthermore, (\ref{cdstar}) implies that $D^* f =  C f$ due to the uniqueness of $D^*$.  
 
(iii) 
Assume that $f \in \dom(S)$. Then $f\in\dom(C)$ 
(see Lemma \ref{dom1}(i.1)). 
Furthermore, $Cf \in \dom(D)$. Therefore, $\dom(S)\subseteq \dom (D  C)$. The converse implication $\dom (D C) \subseteq \dom(S)$ is obvious.
Now it is clear that $S=D  C$.

(iv) Let $c=(c_k)_{k=1}^\infty\in\dom(CD)$. 
Then $\sum_{l=1}^\infty  c_l \<\psi_l,\psi_k\>$ converges for every $k$ and it is equal to $\<Dc,\psi_k\>$. Furthermore, $Dc\in\dom(C)$, which implies that $c\in\dom(G)$ and $Gc=CDc$. 

(v) Since the statement is given in \cite{casoleli1} without proof, for the sake of completeness we refer to \cite[Lemma 3.1]{dslpfr} for a proof. Note that \cite[Lemma 3.1]{dslpfr} concerns sequences satisfying the lower $p$-frame condition, but the proof 
that $C$ is closed does not use the validity of the lower $p$-frame condition. 

(vi) Let $f,g \in \dom(S)$. Then 
$ \left< S f, g\right> = \sum_{k=1}^\infty \left< f , \psi_k \right> \left< \psi_k , g \right> = \left< f, S g \right>,$ which means that $S$ is symmetric. An example of a sequence $\Psi$ with $S$ being non-densely defined can be seen in Remark \ref{rem32}(1). 
Further,
$\left< Sf , f \right> = 
\sum_{k=1}^\infty \left| \left< f , \psi_k\right> \right|^2 \ge 0$ for every $f\in\dom(S)$. 
If $\sum_{k=1}^\infty \left| \left< f , \psi_k\right> \right|^2 = 0$, then 
$f \in (\span\{ \Psi\})^\bot$, which implies that $S$ is positive definite on $\overline{\span} \{ \Psi\} \cap \dom(S)$.

(vii) By (v), $C$ is closed. Assume that $C$ is densely defined. 
Then Proposition \ref{sec:closeunb1}(e) implies that $C^*C$ is self-adjoint and thus, closed.
By (ii) and (iii), we have that $S=DC\subseteq C^*C$, which implies that $S$ is closable.

(viii) and (ix) follow from (i), Lemma \ref{dom1}(iii) and Prop. \ref{sec:closeunb1}(b),(c) applied with $F=D$.

(x) Let $D$ be closed. By Lemma \ref{dom1}(iii) and Proposition \ref{sec:closeunb1}(e), 
it follows that $D^*D$ is self-adjoint, in particular, densely defined. Using (i) and (iv), it follows that $D^*D\subseteq G$, which implies that $G$ is also densely defined.

(xi) Assume that $\sum_{l=1}^\infty |\<\psi_l,\psi_k\>|^2 <\infty$ for every $k\in\mn$. 
Denote $c^n=(c^n_k)_{k=1}^\infty, c=(c_k)_{k=1}^\infty, d=(d_k)_{k=1}^\infty$. Let $c^n\in\dom(G)$, $n\in\mn$,  $c^n\to c$ in $\ell^2$ and $Gc^n\to d$ in $\ell^2$ as $n\to\infty$. 
Fix  $k\in\mn$ and $\varepsilon>0$. Let $N\in\mn$ be such that 
$\|c-c^N\|_{\ell^2} < \frac{\varepsilon}{2} (\sum_l |\<\psi_l,\psi_k\>|^2)^{-1/2}$. 
Since $\sum_{l=1}^\infty c^N_l \<\psi_l,\psi_k\>$ converges, there exists $N_1$ so that $|\sum_{l=i+1}^j c^N_l \<\psi_l,\psi_k\>|< \frac{\varepsilon}{2}$, $\forall i,j >N_1$. Then, for every $i,j >N_1$, 
\begin{eqnarray*} 
|\sum_{l=i+1}^j c_l \<\psi_l,\psi_k\>| &\leq & |\sum_{l=i+1}^j (c_l-c^N_l) \<\psi_l,\psi_k\>| + |\sum_{l=i+1}^j c^N_l \<\psi_l,\psi_k\>| \\
&\leq & \left(\sum_{l=i+1}^j |c_l-c^N_l|^2\right)^{1/2} \left(\sum_{l=i+1}^j |\<\psi_l,\psi_k\>|^2\right)^{1/2} + \frac{\varepsilon}{2} <\varepsilon. 
\end{eqnarray*} 
Therefore, the series $\sum_{l=1}^\infty c_l \<\psi_l,\psi_k\>$ converges.
Furthermore, 
\begin{eqnarray*}
|\sum_{l=1}^\infty (c_l-c^n_l) \<\psi_l,\psi_k\> | &\leq & \|c-c^n\|_{\ell^2} \left(\sum_{l=1}^\infty |\<\psi_l,\psi_k\>|^2\right)^{1/2} \to 0 \mbox{ as } n\to \infty,
\end{eqnarray*}
which implies that $\sum_{l=1}^\infty c^n_l \<\psi_l,\psi_k\> \to \sum_{l=1}^\infty c_l \<\psi_l,\psi_k\>$ as $n\to\infty$. Now, 
since $Gc^n \to d$ in $\ell^2$ as $n\to \infty$ and since convergence in $\ell^2$ implies convergence by coordinates, it follows that $\sum_{l=1}^\infty c^n_l \<\psi_l,\psi_k\> \to d_k$ as $n\to\infty$. 
Therefore,  $\sum_{l=1}^\infty c_l \<\psi_l,\psi_k\>=d_k$ for every $k\in\mn$, 
which implies that $c\in\dom(G)$ and $Gc=d$.
\end{proof}

\begin{rem.} \label{dcstar} {\rm (concerning Prop. \ref{connectop}(ii)): 
If $C$ is densely defined, in general one should not expect the validity of $D=C^*$. 
 A counterexample can be found in \cite{olepinv} after Corollary 4.7. When $\Psi$ is a Bessel sequence for $\Hil$, then $D=C^*$
 (see, e.g., \cite[Theor. 3.2.3 and Lemma 3.2.1]{ole1}). 
}
\end{rem.}

Note that Proposition \ref{connectop}(ix) extends \cite[Corollary 4.7]{olepinv}.

\subsection{Kernels}

\begin{lem.} \label{lemker}
For an arbitary sequence $\Psi$, the following statements hold. 
\begin{itemize}
\item[{\rm (i)}] $\ker(S)=\ker(C)=(\overline{\span}{\{\Psi\}})^{\, \bot}= 
{(\ran(D)) }^{\, \bot}$. 
\item[{\rm (ii)}] $\ker(D)\subseteq \overline{\ran(C)}^{\, \bot}$. 
\end{itemize}
\end{lem.}
\begin{proof} (i) Assume that $f\in\ker(S)$. Then $ \sum_{k=1}^\infty \left| \left< f , \psi_k \right> \right|^2  = \left< S f , f \right>  = 0$, which implies that $f\in\ker(C)$.  It is obvious that $\ker(C)\subseteq\ker(S)$. 
The equality $\ker(C)=\overline{\span}{\{\Psi\}}^{\, \bot}$ is also obvious.
The equality $\ker(C)={\left( \ran(D) \right)}^\bot$ follows from Lemma \ref{dom1}(iii), Proposition \ref{connectop}(i) and Proposition \ref{sec:closeunb1}(d) with $F=D$.

(ii) Follows easily using Prop. \ref{connectop}(i).
\end{proof}

\subsection{Ranges}

\begin{lem.} \label{lemran} For an arbitrary sequence $\Psi$, the following statements hold.
\begin{itemize}
 
 \item[{\rm (i)}]  $\ran(S) \subseteq \ran(D)$. 
\item[{\rm (ii)}] If $\ran(S) = \ran(D)$, then  $\overline{\ran(C)}\cap \dom(D)= \ran(C)\cap\dom(D)$. 
\item[{\rm (iii.1)}] 
Let $\Psi$ be a Bessel sequence for $\Hil$. Then $\ran(S) \subseteq \ran(D)\subseteq \overline{\ran(S)}$ (with dense inclusion). 

 \item[{\rm (iii.2)}]  Let $\Psi$ be a Bessel sequence for $\Hil$. Then $\ran(S) = \ran(D)$ if and only if $\Psi$ is a frame sequence.
\end{itemize} 
  \end{lem.}
    
\begin{proof} (i) follows from Prop. \ref{connectop}(iii). 

(ii)
Assume that $\ran(S)=\ran(D)$. Let $c\in\overline{\ran(C)}\cap \dom(D)$. 
Since $Dc\in\ran(D)=\ran(S)$, there exists $f\in\dom(S)\subseteq \dom(C)$ so that $Dc= DCf$. 
Since $D_{\overline{\ran(C)}\cap \dom(D)}$ is injective (see Lemma \ref{lemker}(ii)), it follows that $c=Cf\in\ran(C)$. Therefore, $\overline{\ran(C)}\cap \dom(D)\subseteq \ran(C)\cap\dom(D)$. The inverse inclusion is obvious.

(iii.1) Let $f\in\ran(D)$. 
Since $\Psi$ is a Bessel sequence for $\Hil$, it follows that $ \ker D=\ker (C^*)=\overline{\ran(C)}^{\,\bot}$, which implies that $f=Dc$ for some $c\in\overline{\ran(C)}$. Therefore, $f=\lim_{n\to\infty} DCf_n=\lim_{n\to\infty} Sf_n$ for some $f_n\in\Hil, n\in\mn$,  such that $c=\lim_{n\to\infty} Cf_n$. 
This completes the proof. 

(iii.2) Assume that $\ran(S)=\ran(D)$. Since $\Psi$ is a Bessel sequence for $\Hil$, it follows that $\dom(D)=\Hil$ 
 and (ii) implies that $\ran(C)$ is closed. Now  
 \cite[Cor. 5.5.3]{ole1} 
 (see Prop.  \ref{sec:framclass2}(b)) implies that $\Psi$ is a frame sequence. Conversely, if $\Psi$ is a frame sequence, then $\ran(S)=\ran(D)=\overline{\span}\{\Psi\}$
 (see, e.g., \cite[Lemma 5.1.5 and Theor. 5.5.1]{ole1}). 
 \end{proof}

 \section{Classification of General Sequences} \label{sec:framclass0}

In this section, we turn to the central topic of the paper, namely, 
the classification of arbitrary sequences. We 
consider two different methods of classification, first in terms of the 
associated operators, then via orthonormal bases.

Note that throughout this section, we use a rather unconventional numbering technique to increase the comparability for the convenience of the reader. 
 The items (a), (b), \ldots \ refer always to the same class of sequences, Bessel sequences, frame sequences, \ldots. Therefore, if no result is available for a given class, we simply omit the corresponding item.

\subsection{Classification by the associated operators}

In the next four propositions, we proceed with the classification of arbitrary sequences in terms of the operators $C, D, S$, and $G$.

\begin{prop.} \label{sec:framclass2} 
Given a sequence $\Psi$, the following statements hold.
\begin{itemize} 
\item[{\rm (a1)}] {\rm \cite[Lemma 11.8]{he98-1}} $\Psi$ is a Bessel sequence for $\Hil$ if and only if $\dom(C) = \Hil$.
\item[{\rm (a2)}]  $\Psi$ is a Bessel sequence for $\Hil$ with bound $B$ if and only if $\dom(C) = \Hil$
and $C$ is bounded with $\|C\| \leq \sqrt{B}$.

\item[{\rm (b)}] 
{\rm\cite[Cor. 5.5.3]{ole1}}
$\Psi$ is a frame sequence if and only if $\dom(C) = \Hil$ and $\ran(C)$ is closed.
\item[{\rm (c)}] 
{\rm\cite[Cor. 5.5.3]{ole1}} $\Psi$ is a frame for $\Hil$ if and only if $\dom(C) = \Hil$, $\ran(C)$ is closed and $C$ is injective.
\item[{\rm (d)}] 
$\Psi$ is a Riesz basis for $\Hil$ if and only if $\dom(C) = \Hil$ and $C$ is bijective. 

\item[{\rm (e)}] {\rm \cite[Lemma 3.1]{casoleli1}} $\Psi$ is a lower frame sequence for $\Hil$ if and only if $C$ is injective and $\ran(C)$ is closed. 

\item[{\rm (f)}] 
{\rm \cite[Ch.4\,Sec.2]{young1}} 
$\Psi$ is a Riesz-Fischer sequence if and only if $C$ is surjective. 
\item[{\rm (g)}] $\Psi$ is complete in $\Hil$ if and only if $C$ is injective.

\end{itemize}
\end{prop.}
\begin{proof} (a2) Having (a1) proved, then (a2) is obvious. 

(d) By \cite[Theor. 5.4.1]{ole1} and  \cite[Prop. 5.1.5]{gr01}, $\Psi$ is a Riesz basis for $\Hil$ if and only if $\Psi$ is a frame for $\Hil$ and $C$ is surjective. The rest follows from (c).

(g) follows from Lemma \ref{lemker}(i). 
\end{proof}

\begin{prop.} \label{sec:classwsynth1} 
Given a sequence $\Psi$, the following statements hold.
\begin{itemize} 
\item[{\rm (a1)}] 
{\rm\cite[Cor. 3.2.4 and Theor. 3.2.3]{ole1}}
$\Psi$ is a Bessel sequence if and only if $\dom(D) = \ell^2$.
\item[{\rm (a2)}] 
{\rm\cite[Theor. 3.2.3]{ole1}} 
$\Psi$ is a Bessel sequence with bound $B$ if and only if $\dom(D) = \ell^2$ and $D$ is bounded with $\|D\|\leq \sqrt{B}$. 
\item[{\rm (b1)}]  
 {\rm\cite[Cor. 5.5.2]{ole1}}
  $\Psi$ is a frame sequence if and only if $\dom(D) = \ell^2$ and $\ran(D)$ is closed.
 \item[{\rm (b2)}] $\Psi$ is a frame sequence if and only if $\ran(D)$ is closed and $\ran(D)\subseteq \dom(C)$.

  \item[{\rm (b3)}] $\Psi$ is a frame sequence if and only if $\dom(D) = \ell^2$ and $\ran(D)=\ran(S)$.

\item[{\rm (c)}]  
   {\rm\cite[Theor. 5.5.1]{ole1}}
  $\Psi$ is a frame if and only if $\dom(D) = \ell^2$ and $D$ is surjective.

  \item[{\rm (d)}]   $\Psi$ is a Riesz basis for $\Hil$ if and only if $\dom(D) = \ell^2$ and $D$ is bijective. 
\item[{\rm (e)}]  
$\Psi$ is a lower frame sequence for $\Hil$ if and only if $ran(D)$ is dense in $\Hil$ and $\ran(D^*)$ is closed. 

\item[{\rm (f)}] 
$\Psi$ is a Riesz-Fischer sequence if and only if $D$ is injective and $D^{-1}$ is bounded on $\ran(D)$. 

\item[{\rm (g)}]  $\Psi$ is complete in $\Hil$ if and only if $\ran(D)$ is dense in $\Hil$.

\end{itemize}
\end{prop.}
\begin{proof} 
(b2) Assume that $\ran(D)$ is closed and $\ran(D)\subseteq \dom(C)$. Then $\Hil=\ran(D) \oplus (\ran(D))^\perp\subseteq \dom(C)$, which by Proposition \ref{sec:framclass2}(a1) implies that $\Psi$ is a Bessel sequence for $\Hil$. Now apply (a1) and (b1). 

Conversely, assume that $\Psi$ is a frame sequence. By (b1), $\ran(D)$ is closed, and clearly, $\dom(C)=\Hil\supseteq \ran(D)$.

(b3) follows from Lemma \ref{lemran}(iii.2) and (a1).

(d) By \cite[Theor. 5.4.1 and 6.1.1]{ole1}, $\Psi$ is a Riesz basis for $\Hil$ if and only if $\Psi$ is a frame for $\Hil$ and $D$ is injective. The rest follows from (c).

(e) 
By Propositions \ref{sec:framclass2}(e) 
and \ref{connectop}(i), $\Psi$ is a lower frame sequence for $\Hil$ 
if and only if $D^*$ is injective and $\ran(D^*)$ is closed. 
Since $D$ is densely defined, Proposition \ref{sec:closeunb1}(d) completes the proof.

(f) is clear (see, e.g., \cite[V.4.3 Theorem 2 and the Remark after that]{ka84}). 

(g) follows from Lemma \ref{lemker}(i).
 \end{proof}

Note that $(f)$ (a result about the Riesz-Fischer sequences) is missing in Proposition \ref{sec:framclass3}. 
The operator S does not distinguish between a frame which is a Riesz-Fischer sequence and a frame which is not a Riesz-Fischer sequence. For example, an orthonormal basis and a Parseval frame have the same frame operator, the identity. This can also be seen in Proposition \ref{sec:framclass3}(d), where the properties of $(S^{-1}\psi_k)$ have to be taken into account.

\begin{prop.} \label{sec:framclass3}
Given a sequence $\Psi$, the following statements hold. 
\begin{itemize} 
\item[{\rm (a1)}] $\Psi$ is a Bessel sequence for $\Hil$ if and only if $\dom(S)=\Hil$. 
\item[{\rm (a2)}]  $\Psi$ is a Bessel sequence for $\Hil$ with bound $B$ if and only if $\dom(S)=\Hil$ and 
$S$ is bounded with $\|S\| \leq B$.

\item[{\rm (b1)}]  $\Psi$ is a frame sequence if and only if $\dom(S)=\Hil$ and $\ran(S)$ is closed. 

\item[{\rm (b2)}]   $\Psi$ is a frame sequence if and only if $\dom(S)=\Hil$ and $\ran(S)=\overline{\span}\{\Psi\}$. 

\item[{\rm (c1)}] {\rm \cite[Theorem 2.1]{ch93-2}}  $\Psi$ is a frame for $\Hil$ if and only if $\dom(S) = \Hil$ and $S$ is surjective.

\item[{\rm (c2)}]  
$\Psi$ is a frame for $\Hil$ if and only if $\dom(S) = \Hil$ and $S$ is  bijective.

\item[{\rm (d)}]  
$\Psi$ is a Riesz basis for $\Hil$ if and only if $\dom(S)=\Hil$, $S$ is bijective and $(S^{-1}\psi_k)_{k=1}^\infty$ is biorthogonal to $\Psi$.

\item[{\rm (e)}] $\Psi$ is a lower frame sequence for $\Hil$ if and only if $S$ is injective and $\ran(C)$ is closed. 

\item[{\rm (g)}]  $\Psi$ is complete in $\Hil$ if and only if $S$ is injective. 
\end{itemize}
\end{prop.}
\begin{proof} (a1) 
Let $\dom(S)=\Hil$. Validity of the upper frame inequality follows from the first part of the proof in \cite[Theorem 2.1]{ch93-2}. 
For the sake of completeness, we sketch a proof. If $\dom(S) = \Hil$, then $\sum|\<f,\psi_k\>|^2=\<Sf,f\> <\infty$ for every $f\in\Hil$, which implies that $\dom(C)=\Hil$ and the rest follows from Proposition \ref{sec:framclass2}(a1). 
For the other direction, see 
the text after Example 5.1.4 in \cite[Sec.\,5.1]{ole1}.

(a2) 
See (a1); the boundedness of $S$ can be found, e.g., in the proof of \cite[Lemma 5.1.5]{ole1}. 

(b1) 
Assume that $\dom(S)=\Hil$ and $\ran(S)$ is closed. By (a1), $\Psi$ is a Bessel sequence for $\Hil$. Then Lemma \ref{lemran}(iii.1) implies that $\ran(D)$ is closed, which by Proposition \ref{sec:classwsynth1}(b1) implies that $\Psi$ is a frame sequence. 

Conversely, if $\Psi$ is a frame sequence, 
the conclusion follows from \cite[Theorem 2.1]{ch93-2} (see (c1)).

(b2) follows from \cite[Theorem 2.1]{ch93-2} (see (c1)).

(c2) follows from (c1) and \cite[Lemma 5.1.5]{ole1}. 

(d) Use (c2) and the fact that $\Psi$ is a Riesz basis for $\Hil$ if and only if $\Psi$ is a frame for $\Hil$ and $(S^{-1}\psi_k)_{k=1}^\infty$ is biorthogonal to $\Psi$ \cite[Theor. 5.4.1 and 6.1.1]{ole1}.

(e) By Lemma \ref{lemker}(i), $S$ is injective if and only if $C$ is injective. The rest follows from Proposition \ref{sec:framclass2}(e).

(g)  follows from Lemma \ref{lemker}(i).
\end{proof}

\begin{prop.} \label{sec:classgram1} 
For a sequence $\Psi$, the following statements hold. 
\begin{itemize}
\item[{\rm(a1)}] $\Psi$ is a Bessel sequence for $\Hil$ if and only if $\dom(G)=\ell^2$. 
\item[{\rm(a2)}] {\rm \cite[Lemma 3.5.1]{ole1}}  $\Psi$ is a Bessel sequence for $\Hil$ with bound $B$ if and only if $\dom(G)=\ell^2$ and $G$ is bounded with $\|G\|\leq B$.

\item[{\rm(a3)}]
If $\Psi$ is a Bessel sequence for $\Hil$, then $G\vert_{\ran(C)}$ is an injective operator from $\ran(C)$ into $\ran(C)$ and $\ran(G\vert_{\ran(C)})$ is dense in $\ran(C)$. 

\item[{\rm(b)}] 
$\Psi$ is a frame sequence if and only if $\dom(G)=\ell^2$ and $G\vert_{\ran(C)}$ is a bounded operator from $\ran(C)$ onto $\ran(C)$ with bounded inverse.

\item[{\rm(c)}] 
$\Psi$ is a frame for $\Hil$ if and only if $\Psi$ is complete in $\Hil$, $\dom(G)=\ell^2$ and $G\vert_{\ran(C)}$ is a bounded operator from $\ran(C)$ onto $\ran(C)$ with bounded inverse.

\item[{\rm(d)}] 
{\rm \cite[Theor. 3.6.6]{ole1}} $\Psi$ is a Riesz basis for $\Hil$ if and only if $\Psi$ is complete in $\Hil$ and $G$ is a bounded invertible operator on $\Hil$.

\item[{\rm(f)}] {\rm \cite[Ch.4 Sec.2]{young1}} $\Psi$ is a Riesz-Fischer sequence
with bound $A$ if and only if the \lq\lq sections\rq\rq\ $G_n$ of $G$ 
satisfy the inequality  $A \norm{}{c} \le \norm{}{ G_n c}$ for all finite sequences and all $n$, where the \lq\lq section\rq\rq\ $G_n$ is the matrix $(\<\psi_l, \psi_k\>)_{l,k=1}^n$.

\end{itemize}
\end{prop.}
\begin{proof} 
(a1) 
Let $\dom(G)=\ell^2$. Then $G$ admits a matrix representation via the orthonormal basis $(\delta_k)_{k=1}^\infty$ of $\ell^2$ and thus, $G$ is bounded \cite[Sec. 29 Theor. 1]{AhiezerGlazman}.\footnote{For a matrix representation using frames refer to \cite{xxlframoper1}.}
Since \cite{AhiezerGlazman} is in Russian, for convenience of the reader we add a sketch of a proof.  
Let $n\in\mn$. Consider the operator $G_n$ defined by $G_n(c)_{c=1}^\infty=(\sum_{l=1}^\infty c_l\<\psi_l,\psi_k\>)_{k=1}^n$. Then $\dom(G_n)=\ell^2$ and $G_n$ is bounded. 
Therefore, by the Banach-Steinhaus theorem, $G$ is bounded. 
Now the rest follows from {\rm \cite[Lemma 3.5.1]{ole1}}  (see (a2)).

(a3) follows from \cite[Lemma 3.5.2]{ole1} and Proposition \ref{connectop}(i).

(b) If $\Psi$ is a frame sequence, the statement 
follows from \cite[Prop.\,5.2.2]{ole1} and Proposition \ref{connectop}(i).
Conversely, assume that $\dom(G)=\ell^2$ and $G\vert_{\ran(C)}$ is a bounded operator from $\ran(C)$ onto $\ran(C)$ with bounded inverse.
Then $\Psi$ is a Bessel sequence for $\Hil$ and thus $G=CD$.
Consider the subspace $\overline{\span}\{\Psi\}$ of $\Hil$ and 
the operator $C_1=C_{\overline{\span}\{\Psi\}}$. 
By Lemma \ref{lemker}(i), 
$C_1$ is a bijection from $\overline{\span}\{\Psi\}$ onto $\ran(C)$.
Since $\ran(D)\subseteq \overline{\span}\{\Psi\}$, 
for every $c\in\ran(C)$ one has $G\vert_{\ran(C)} c=CDc=C_1D\vert_{\ran(C)}c$ and thus, 
$C_1 D\vert_{\ran(C)} (G\vert_{\ran(C)})^{-1}d =d$ for every $d\in\ran(C)$. 
Therefore, $$C_1^{-1}=D\vert_{\ran(C)} (G\vert_{\ran(C)})^{-1},$$
which implies that $C_1^{-1}$ is bounded. 
Thus, the bounded injective operator $C_1$ has a bounded inverse on $\ran(C_1)$, which implies that $C_1$ has closed range \cite[Cor. IV.3.6]{wern1} (for a reference in English, see \cite[Ex. VI.9.15(i)]{DunfordS1}). 
 Therefore, by Proposition \ref{sec:framclass2}(c), $\Psi$ is a frame for $\overline{\span}\{\Psi\}$.
 
 (c) follows from (b).
\end{proof}

Note that 
again some items are missing, as in Proposition \ref{sec:framclass3}. It is impossible to have some connection of properties of the Gram matrix to completeness. For example, the standard orthonormal basis $(e_1,e_2,e_3,\ldots)$ as well as $(e_2,e_3,e_4,....)$ have the same Gram matrix. 
The problem concerning sequences that satisfy the lower frame condition is still open.

  \subsection{Operators that preserve the sequence properties}\label{sec:oponfram}  

\begin{prop.} \label{classframes} 
Given a sequence $\Psi$, the following statements hold. 
\begin{itemize} 
\item[{\rm (a)}] \label{sec:besselop1} 
Let $\Psi$ be a Bessel sequence for $\Hil$ with bound $B$ and let $F : \span\{\Psi\} \to \Hil$ be a bounded operator. Then $(F \psi_k)_{k=1}^\infty$ is a Bessel sequence for $\Hil$ with bound  $B \, \|F\|^2$.

\item[{\rm (b)}] \label{sec:framseqop1}  
Let $\Psi$ be a frame sequence 
with bounds $A,B$ and let $F : \overline{\span}\{\Psi\}\to\overline{\span}\{\Psi\}$ be a bounded surjective operator. Then $(F \psi_k)_{k=1}^\infty$ is a frame sequence with frame bounds $A  \|{F^\dagger}\|^{-2}$ and $B  \|F\|^{2}$, where $F^\dagger$ \footnote{For the pseudo-inverse of a bounded operator with closed range see, e.g., \cite[App. A.7]{ole1}.} denotes the pseudo-inverse of $F$.

\item[{\rm (c)}] \label{sec:frameop1} \index{frames!operators applied on}
{\rm \cite[Cor. 5.3.2]{ole1} } Let $\Psi$ be a frame for $\Hil$ with bounds $A,B$ and let $F : \Hil\to\Hil$ be a bounded surjective operator. Then $(F \psi_k)_{k=1}^\infty$ is a frame for $\Hil$ 
with bounds $A  \|{F^\dagger}\|^{-2}$ and $B  \|F\|^{2}$, where $F^\dagger$ denotes the pseudo-inverse of $F$.

\item[{\rm (d)}] Let $\Psi$ be a Riesz basis for $\Hil$ with bounds $A,B$ and let $F : \Hil\to\Hil$ be a bounded bijective operator. Then  $(F \psi_k)_{k=1}^\infty$ is a Riesz basis for $\Hil$ with bounds $A\|F^{-1}\|^{-2}$, $B\|F\|^2$.

\item[{\rm (e)}] Let $\Psi$ be a lower frame sequence for $\Hil$ with bound $A$ and let the operator $F:\overline{\span\{\Psi\}}\to\Hil$ be bounded and such that $F^*$ has a bounded inverse on $\ran(F^*)$ (equivalently, $F$ be bounded and surjective).  Then  $(F \psi_k)_{k=1}^\infty$ is a lower frame sequence for $\Hil$ with bound $A\| (F^*)^{-1}\|^{-2}$.

\item[{\rm (f)}] Let $\Psi$ be a Riesz-Fischer sequence with bound $A$ and let 
 the operator $F:\span\{\Psi\}\to\Hil$ have a bounded inverse $F^{-1}:\ran(F)\to\Hil$ with bound $K$. 
 Then  $(F \psi_k)_{k=1}^\infty$ is a Riesz-Fischer sequence with bound $AK^{-2}$.

\item[{\rm (g)}] Let $\Psi$ be complete in $\Hil$ and let $F : \span\{\Psi\}\to\Hil$ be bounded with $\ran(F)$ being dense in $\Hil$. Then $(F\psi_k)_{k=1}^\infty$ is complete in $\Hil$.

\end{itemize}
\end{prop.}
\begin{proof} (a) 
Let $\overline{F}$ denote the bounded extension of $F$ on $\overline{\span}\{\Psi\}$ without increasing the norm. For every $f\in\Hil$, $$ \sum_{k=1}^\infty \left| \left< f, F \psi_k \right> \right|^2 = \sum_{k=1}^\infty \left| \left< \overline{F}^{\,*} f, \psi_k \right> \right|^2 \le B \, \|\overline{F}^{\,*} f\|^2 \le B \, \|F\|^2\cdot\|f\|^2. $$

(b) follows from \cite[Cor. 5.3.2]{ole1} (see (c)).

(d) 
By (c), $(F\psi_k)_{k=1}^\infty$ is a frame for $\Hil$ and thus, complete in $\Hil$ (see Proposition \ref{sec:classwsynth1}(c)). Furthermore, 
for every finite sequence $(c_k)$, 
$$ \|\sum c_k \, F\psi_k\|^2= \|F(\sum c_k \, \psi_k)\|^2 \geq \frac{1}{\|F^{-1}\|^2} A\sum |c_k|^2$$
and similar, $\|\sum c_k \, F\psi_k\|^2 \leq \|F\|^2 \,B\sum |c_k|^2$.

(e) 
A bounded operator from a Banach space into a Banach space is surjective if and only if its adjoint has a bounded inverse on the range of the adjoint (see, e.g., \cite[Theorem 4.15]{rudinengl} and \cite[V.4.3 Theorem 2 and the Remark after that]{ka84}). 
 For every $f\in\Hil$, 
$$ \sum_{k=1}^\infty \left| \left< f, F \psi_k \right> \right|^2 = \sum_{k=1}^\infty \left| \left< F^* f, \psi_k \right> \right|^2 \geq A \, \|F^* f\|^2 \geq \frac{A}{\| (F^*)^{-1}\|^2}\,\|f\|^2.$$

(f) As in (d),  $ \|\sum c_k \, F\psi_k\|^2 \geq K^{-2}A\sum |c_k|^2$ for every finite sequence $(c_k)$.

(g) Let $\Psi$ be complete in $\Hil$ and let $\overline{F}$ denote the bounded extension of $F$ on $\overline{\span}\{\Psi\}=\Hil$. 
 Assume that $f\in\Hil$ and $\<f,F\psi_k\>=0, \forall k\in\mn$. Then $\overline{F}^{\,*} f=0$. Since $\ran(\overline{F})$ is dense in $\Hil$, it follows that $\overline{F}^{\,*}$ is injective 
 (see Proposition \ref{sec:closeunb1}(d)). Therefore, $f=0$.
\end{proof}

\subsection{Classification with Orthonormal Bases} \label{sec:classwonbs0}

Another way of classifying the sequences is to examine how an orthonormal basis behaves under application of a given class of operators.
This we do in the present subsection. From now on, $D_\Psi$ denotes the synthesis operator for $\Psi$ and $C_{(e_k)}$ denotes the analysis operator for $(e_k)_{k=1}^\infty$.

\begin{prop.} \label{sec:seqonb1} Let $(e_k)_{k=1}^{\infty}$ be an orthonormal basis for $\Hil$. 
\begin{itemize} 
\item[{\rm (a)}]  
The Bessel sequences for $\Hil$ are precisely the families $(V e_k)_{k=1}^\infty$, where $V : \Hil \rightarrow \Hil$ is a bounded operator. 
\item[{\rm (b)}] 
The frame sequences for $\Hil$ are precisely the families $(V e_k)_{k=1}^\infty$, where $V : \Hil \rightarrow \Hil$ is a bounded operator with closed range.
\item[{\rm (c)}] \index{frame sequences!classification with ONBs} 
{\rm \cite[Theor. 5.5.5]{ole1}} 
The frames for $\Hil$ are precisely the families $(V e_k)_{k=1}^\infty$, where $V : \Hil \rightarrow \Hil$ is a bounded and surjective operator. 
\item[{\rm (d)}] {\rm \cite[Def. 3.6.1 and Theor. 3.6.6]{ole1}} The Riesz bases for $\Hil$ are precisely the sequences $(V e_k)_{k=1}^\infty$, where $V:\Hil\to\Hil$ is a bounded bijective operator.

\item[{\rm (e)}] 
The lower frame sequences for $\Hil$ are precisely the families $(Ve_k)_{k=1}^\infty$, 
where $V : \dom(V) \to \Hil$ is a densely defined operator such that $e_k\in\dom(V)$, $\forall k\in\mn$, $V^*$ is injective with bounded inverse on $\ran(V^*)$, and $V(\sum_{k=1}^n c_k e_k)\to V(\sum_{k=1}^\infty c_k e_k)$ as $n\to\infty$ for every $\sum_{k=1}^\infty c_k e_k\in\dom(V)$.

\item[{\rm (f)}] 
 {\rm \cite[Prop. 2.3]{casoleli1}}  The Riesz-Fischer sequences 
 are precisely the families $(V e_k)_{k=1}^\infty$, where $V$ is an operator having all $e_k$ in the domain and which has a bounded inverse $V^{-1}:\ran(V)\to\Hil$. 

\item[{\rm (g)}] 
The complete sequences are precisely the families $(Ve_k)_{k=1}^\infty$, 
where $V : \dom(V) \to \Hil$ is a densely defined operator such that $e_k\in\dom(V)$, $\forall k\in\mn$, $\ran(V)$ is dense in $\Hil$ (equivalently, the adjoint $V^*$ is injective) and $V(\sum_{k=1}^n c_k e_k)\to V(\sum_{k=1}^\infty c_k e_k)$ as $n\to\infty$ for every $\sum_{k=1}^\infty c_k e_k\in\dom(V)$.

\end{itemize}
\end{prop.} 
\begin{proof} (a) 
If $V:\Hil\to\Hil$ is a bounded operator, the conclusion follows from Proposition \ref{classframes}(a). 

Conversely, assume that $\Psi$ is a Bessel sequence for $\Hil$. 
Then the operator $V = D_{\Psi}  C_{(e_k)}$ is bounded and $\Psi=(Ve_k)_{k=1}^\infty$.

(b) If $V : \Hil \rightarrow \Hil$ is a bounded operator with closed range, then \cite[Prop. 5.3.1]{ole1} implies that $(Ve_k)_{k=1}^\infty$ is a frame sequence.

Conversely,  assume that $\Psi$ is a frame sequence.  
Consider the operator $V = D_{\Psi} C_{(e_k)}$. Since $C_{(e_k)}$ is a bounded bijection (see Proposition \ref{sec:framclass2}(d)) and $D_{\Psi}$ is a bounded operator with closed range (see Proposition \ref{sec:classwsynth1}(b1)), it follows that $V$ is bounded and $\ran(V)$ is closed. Clearly, $\Psi=(Ve_k)_{k=1}^\infty$.

(e) Let $\Psi$ be a lower frame sequence for $\Hil$. Consider the operator $V:=D_\Psi C_{(e_k)}$. It is clear that $\dom(V)$ contains all $e_k$, $k\in\mn$, (and thus, $V$ is densely defined) and $(Ve_k)_{k=1}^\infty=\Psi$. 
Using the fact that $C_{(e_k)}$ is a bounded bijection of $\Hil$ onto $\ell^2$, and $C_{(e_k)}$ maps bijectively $\dom(V)$ onto $\dom(D)$, it is not difficult to see that $\dom (V^*)= \dom ((D_\Psi)^*) =\dom ((C_{(e_k)})^* (D_\Psi)^*)$ 
and $V^* = (C_{(e_k)})^* (D_\Psi)^*$. Using Proposition \ref{connectop}(i), it follows that $V^*=D_{(e_k)} C_\Psi$. Now Propositions \ref{sec:framclass2}(e) and \ref{sec:classwsynth1}(d)  imply that $V^*$ is injective. 
Furthermore, $(V^*)^{-1}=(C_\Psi)^{-1} (D_{(e_k)})^{-1}\vert_{\ran(D_{(e_k)} C_\Psi)}$
 and $(V^*)^{-1}$ is bounded on $\ran(V^*)$.
Now assume that $\sum_{i=1}^\infty c_i e_i\in\dom(V)$. Then $(c_k)_{k=1}^\infty=C_{(e_k)}(\sum_{i=1}^\infty c_i e_i)\in\dom(D_\Psi)$, which implies that $\sum_{k=1}^\infty c_k\psi_k$ converges in $\Hil$ and $V(\sum_{i=1}^\infty c_i e_i)=\sum_{k=1}^\infty c_k\psi_k$. For every $n\in\mn$, 
$\sum_{k=1}^n c_k e_k\in\dom(V)$ and $V(\sum_{k=1}^n c_k e_k)=\sum_{k=1}^n c_k \psi_k$. 
Therefore, $V(\sum_{k=1}^n c_k e_k)\to V(\sum_{i=1}^\infty c_i e_i)$ as $n\to\infty$.

Conversely, assume that $V : \dom(V) \to \Hil$ satisfies the conditions in the statement of (e). 
Let $f\in\Hil$ be such that $\sum_{k=1}^\infty |\<f,Ve_k\>|^2<\infty$ and denote $N:=\sum_{k=1}^\infty |\<f,Ve_k\>|^2$. First we prove that $f$ belongs to $\dom(V^*)$. 
Let $\sum_{k=1}^\infty c_k e_k$ be an arbitrary element of $\dom(V)$. 
For every $n\in\mn$,  $\sum_{k=1}^n c_k e_k$ belongs to $\dom(V)$ and
\begin{eqnarray*}
\left|  \<V \left(\sum_{k=1}^n c_k e_k\right) , f\> \right| 
&\leq& \left(\sum_{k=1}^n |c_k|^2\right)^{1/2} \left(\sum_{k=1}^n |\<Ve_k,f\>|^2\right)^{1/2} \\
&\leq& \sqrt{N} \left\|\sum_{k=1}^n c_k e_k \right\| . 
\end{eqnarray*}
By the assumptions, $V(\sum_{k=1}^n c_k e_k)\to V(\sum_{k=1}^\infty c_k e_k)$ as $n\to\infty$. Therefore, taking the limit as $n\to\infty$ in the above inequalities, it follows that
\begin{equation*}
\left|  \<V \left(\sum_{k=1}^\infty c_k e_k\right) , f\> \right| \leq \sqrt{N} \left\|\sum_{k=1}^\infty c_k e_k \right\|. 
\end{equation*}
Hence, $f\in\dom(V^*)$ and 
\begin{equation*}\label{lowerf}
\sum_{k=1}^\infty |\<f,Ve_k\>|^2=\|V^*f\|^2\geq \frac{1}{\|(F^*)^{-1}\|} \|f\|^2.
\end{equation*}

(f) Since the statement in {\rm \cite[Prop. 2.3]{casoleli1}} is given without proof, for the sake of completeness we add a proof here. One of the directions follows from Proposition  \ref{classframes}(f). 
For the other direction, assume that $\Psi$ is a Riesz-Fischer sequence 
with bound $A$.  
 Consider then the operator $V = D_{\Psi} C_{(e_k)}$. Then the domain of $V$ contains 
 all $e_k$, 
 $(Ve_k)_{k=1}^\infty=\Psi$, 
 and $\|Vf\|=\|D_{\Psi} C_{(e_k)}f\|\geq \sqrt{A}\|C_{(e_k)}f\|=\sqrt{A}\|f\|$  for every $f\in\dom(V)$.  Therefore, $V$ is injective with bounded inverse on $\ran(V)$.

 (g) Let $\Psi$ be complete in $\Hil$. Consider the operator $V:=D_\Psi C_{(e_k)}$. In a similar way as in (e), it follows that $V$ has the desired properties.

Conversely, assume that $V : \dom(V) \to \Hil$ satisfies the conditions in the statement of (g). 
Let $f\in\Hil$ be such that $\<f,Ve_k\>=0$, $\forall k\in\mn$. First we prove that $f$ belongs to $\dom(V^*)$. Every finite sum $\sum c_k e_k$ belongs to $\dom(V)$ and
$ \<V \left(\sum c_k e_k\right) , f\>=0$.
Now, let $\sum_{k=1}^\infty c_k e_k\in\dom(V)$. 
By the above, it follows that 
$  \<V \left(\sum_{k=1}^n c_k e_k\right) , f\> =0, \ \forall n\in\mn$.
By assumption, $V(\sum_{k=1}^n c_k e_k)\to V(\sum_{k=1}^\infty c_k e_k)$ as $n\to\infty$. Therefore, $  \<V \left(\sum_{k=1}^\infty c_k e_k\right) , f\> =0$.
Hence, $   \<V h, f\> =0, \ \forall h\in\dom(V) $, which implies that $f\in\dom(V^*)$.
Therefore, the assumptions $\<f,Ve_k\>=0$, $\forall k\in\mn$, imply that $\<V^*f,e_k\>=0$, $\forall k\in\mn$, which implies that $V^*f=0$ and the injectivity of $V^*$ implies that $f=0$. This ends the proof that $(Ve_k)_{k=1}^\infty$ is complete in $\Hil$.
\end{proof}

Note that, in each of the statements in Proposition \ref{sec:seqonb1}, the orthonormal basis can be replaced by a Riesz basis, in view of Proposition \ref{sec:seqonb1}(d).

\vspace{.25in}
\centerline{ACKNOWLEDGEMENT}

\vspace{.1in}
 The authors would like to thank M. A. El-Gebeily for helpful comments and suggestions.  
This work was partly supported by the WWTF project MULAC (`Frame Multipliers: Theory and Application in Acoustics', MA07-025).
The first two authors are thankful for the hospitality of the 
Institut de Recherche en Math\'ematique et  Physique, Universit\'e catholique de Louvain. 
The second author is grateful for the support from the MULAC-project. She also thanks 
the University of Architecture, Civil Engineering and Geodesy, and the Department of Mathematics of UACEG for supporting the present research.

\end{document}